\documentclass[graybox]{svmult}
\usepackage{mathptmx}       % selects Times Roman as basic font
\usepackage{helvet}         % selects Helvetica as sans-serif font
\usepackage{courier}        % selects Courier as typewriter font
\usepackage{type1cm}        % activate if the above 3 fonts are
\usepackage{makeidx}         % allows index generation
 
\usepackage{multicol}        % used for the two-column index
\usepackage[bottom]{footmisc}% places footnotes at page bottom

\usepackage{amssymb}                 % for mathbb fonts etc
\usepackage{amsmath}
\usepackage{bm,enumerate}
\usepackage[]{graphicx}
\usepackage{epsfig}
\usepackage{subfigure}
\usepackage{cite}
\usepackage{multirow}
\def\BARIT#1{{\bar {#1}}}
\def\BARH{\BARIT{H}}
\def\LATT{{\Lambda}_N}

\def\LATTC{{\bar{\Lambda}_{M}}}

\def\VAR{\mathrm{Var}\,}

\def\Bb{\mathcal{B}}

\def\Rr{\mathcal{R}}

\def\Kk{\mathcal{K}}
\def\Oo{\mathcal{O}}

\def\Qq{\mathcal{Q}}

\def\E{\mathbb{E}}

\def\EXP#1{e^{#1}}

\def\EXPECT{{\mathbb{E}}}

\def\COMMA{\,,}             % use for commas at the ends of formulae
\def\PERIOD{\,.}            % use for periods at the ends of formulae
\def\SEP{{\,|\,}}           % use for seperator in set defs, like {x\in R \SEP

\def\BIGO{\Oo}

\newtheorem{algorithm}{Algorithm}%[section]
\newtheorem{thm}{Theorem}

\newtheorem{defin}{Definition}

\makeindex             % used for the subject index
\begin{document}

%%%%%%%%%%%%%%%%%%%%%%%%%%%%%%%%%%%%%%%%%%%%%%%%%%%%%%%%%%%%%%%%%%%
%
%------------------------------------------------------------------
% Title and authors
%------------------------------------------------------------------
%
\title*{Coupled coarse graining and Markov Chain Monte Carlo for lattice systems}

             \author{Evangelia Kalligiannaki\thanks{The research of E.K. was supported by the National Science Foundation under the grant NSF-CMMI-0835582.} and Markos A. Katsoulakis\thanks{The research of M.A. K. was supported by the National Science Foundation through the grants NSF-DMS-0715125
and the CDI -Type II award NSF-CMMI-0835673.} and Petr Plech\'a\v{c}\thanks{The research of P.P. was partially supported by the National Science Foundation under the grant
NSF-DMS-0813893 and by the Office of Advanced Scientific Computing Research,
U.S. Department of Energy under DE-SC0001340; the work was partly done at the Oak Ridge National Laboratory, 
which is managed by UT-Battelle, LLC under Contract No. DE-AC05-00OR22725.  }}
             \institute{Evangelia Kalligiannaki \at Joint Institute for Computational Sciences, 
        University of Tennessee and Oak Ridge National Laboratory, \email{evy@ornl.gov}
\and Markos A. Katsoulakis \at University of Massachusetts, Department of Applied Mathematics, 
            University of Crete and  Foundation of Research and Technology-Hellas, Greece, \email{ markos@math.umass.edu} \and Petr Plech\'a\v{c} \at Department of Mathematics,
             University of Tennessee, \email{plechac@math.utk.edu} }


\maketitle

%%%%%%%%%%%%%%%%%%%%%%%%%%%%%%%%%%%%%%%%%%%%%%%%%%%%%%%%%%%%%%%%%%%
%
%%%%%%%%%%%%%%%%%%%%%%%%%%%%%%%%%%%%%%%%%%%%%%%%%%%%%%%%%%%%%%%%%%%
%
%  ABSTRACT
%
%%%%%%%%%%%%%%%%%%%%%%%%%%%%%%%%%%%%%%%%%%%%%%%%%%%%%%%%%%%%%%%%%%%

\abstract{
We propose an efficient Markov Chain Monte Carlo method for sampling equilibrium distributions 
for stochastic lattice models, capable of handling correctly long and short-range particle interactions.
The proposed method is a Metropolis-type algorithm with the proposal 
probability transition matrix  based on  the coarse-grained approximating measures introduced in
\cite{KMV1,KPRT}. 
We prove that the  proposed algorithm reduces the computational cost due to energy differences  
and has comparable mixing properties with the classical microscopic Metropolis algorithm, 
controlled by the level of coarsening and reconstruction procedure. 
The properties and effectiveness of the algorithm are demonstrated  with an exactly solvable  example of a one dimensional Ising-type  model, 
comparing  efficiency of the single spin-flip Metropolis dynamics  and the proposed coupled Metropolis algorithm. 
}

%%%%%%%%%%%%%%%%%%%%%%%%%%%%%%%%%%%%%%%%%%%%%%%%%%%%%%%%%%%%%%%%%%%%  
%
% Introduction 
%
%%%%%%%%%%%%%%%%%%%%%%%%%%%%%%%%%%%%%%%%%%%%%%%%%%%%%%%%%%%%%%%%%%%%
%

\section{Introduction}
 \label{intro}
Microscopic, {\em extended\,} (many-particle) systems with {\em complex interactions}   
are ubiquitous in science and engineering applications in a variety of physical and chemical systems,  exhibiting rich mesoscopic  morphologies. For example,  
nano-pattern formation via self-assembly, arises  in surface processes e.g., in  heteroepitaxy,  induced by  competing  short and long-range  interactions \cite{D1}. Other   examples include  macromolecular systems such as polymers, proteins and other   soft matter systems, quantum dots  and
micromagnetic materials. Scientific computing for this class of  systems can rely on molecular simulation methods such as Kinetic Monte Carlo (KMC) or Molecular Dynamics (MD),  however their extensivity, their  inherently complex interactions and   
stochastic nature,  severely limit the spatio-temporal scales that can be addressed by these direct numerical simulation methods.

One of our  primary goals  is to develop  systematic mathematical and computational 
strategies for the speed-up of  microscopic   simulation methods by developing  
{\em coarse-grained} (CG) approximations,  thus reducing the extended system's degrees of freedom.
To date coarse-graining methods have been a subject of intense focus, mainly outside mathematics and primarily in  the physics, applied sciences and engineering literatures \cite{kremerplathe, mp, KMV1, g98, g219}. 
The existing approaches can
give unprecedented speed-up to molecular simulations and can work well
in certain parameter regimes, for instance,  at high temperatures or low density. 
On the other hand,   in many  parameter regimes,  important macroscopic properties may not be  captured properly, e.g. \cite{AKr, mp, PK06}.
Here  we propose to, 
 develop  
{\em reliable} CG algorithms for 
stochastic lattice systems  with complex, and often competing particle interactions  in  equilibrium.  
\noindent
Our proposed methodologies 
stem from the synergy of stochastic processes, 
statistical mechanics and  statistics sampling methods.

%%%%%%%%%%%%%%%%%%%%%%%%%%%%%%%%%%%%%%%%%%%%%%%%%%%%%%%%%%%%%%%%%%%%%%%%%%%%%%%%%%%%%%%%%%%%%
%  Monte Carlo Methods
%%%%%%%%%%%%%%%%%%%%%%%%%%%%%%%%%%%%%%%%%%%%%%%%%%%%%%%%%%%%%%%%%%%%%%%%%%%%%%%%%%%%%%%%%%%%%
Monte Carlo algorithms provide a computational tool capable of estimating observables defined on high-dimensional configuration spaces
that are typical for modeling of complex interacting particle systems  at  or out of equilibrium.
  Markov Chain Monte Carlo (MCMC) simulation methods such as the Metropolis algorithm, were first  proposed in 1953 by Metropolis and his coauthors \cite{MRR} 
for the numerical calculation of the equation of state for a system of rigid spheres. 
It was  generalized  in 1970 by Hastings \cite{HA} and it is commonly referred to as  the Metropolis-Hastings (MH) Monte Carlo method. 
This method belongs to the family of  MCMC methods  which generate ergodic Markovian chains with the
 stationary distribution being the desired sampled probability measure.
Metropolis algorithm consists of two main ingredients: (a) the  probability transition kernel  $q$; the {\it proposal}, 
that generates trial states   and  (b)  the  acceptance probability 
$\alpha$ according  to which the proposed trial is accepted or rejected.
There are though some drawbacks of this method when applied to   large systems, such as a small acceptance probability $\alpha$, 
that leads to costly calculations of a large number of samples that are discarded.
A way to reduce these costs is to {\it predict} efficient  proposal measures such that the computational cost of calculating 
a sample is lower and, if possible, increase the acceptance probability.
%%%%%%%%%%%%%%%%%%%%%%%%%%%%%%%%%%%%%%%%%%%%%%%%%%%%%%%%%%%%%%%%%%%%%%%%%%%%%%%%%%%%%%%%%%%%%
% Mathematics, Convergence, Mixing time
%%%%%%%%%%%%%%%%%%%%%%%%%%%%%%%%%%%%%%%%%%%%%%%%%%%%%%%%%%%%%%%%%%%%%%%%%%%%%%%%%%%%%%%%%%%%%
Convergence and ergodicity properties   of Metropolis type algorithms are studied extensively in a series of works \cite{PDSC1,DSC,RC}. The rate of convergence to stationarity is strongly dependent on the proposal distribution and its relation to the stationary measure (\cite{RC} ch. 7). A quantity that measures the speed of convergence in distribution to stationarity is the spectral gap.  In order to improve an MCMC method one has to increase its spectral gap by smartly constructing a good proposal.
%%%%%%%%%%%%%%%%%%%%%%%%%%%%%%%%%%%%%%%%%%%%%%%%%%%%%%%%%%%%%%%%%%%%%%%%%%%%%%%%%%%%%%%%%%%%%
% Motivate/Introduce our work
%%%%%%%%%%%%%%%%%%%%%%%%%%%%%%%%%%%%%%%%%%%%%%%%%%%%%%%%%%%%%%%%%%%%%%%%%%%%%%%%%%%%%%%%%%%%%

In this work we propose the Coupled Coarse Graining Monte Carlo (Coupled CGMC) method, a new method of constructing  efficient  {\it proposal measures}  based on coarse-graining properties of 
the sampling models.   
We   prove that such approach is  suitable  for models that include both short and long-range interactions between particles.
 Long-range interactions are well-approximated by coarse graining techniques \cite{KMV1, KPR, KT}, and 
Coarse Graining Monte Carlo (CGMC) are adequate simulation methods with an order of
 acceleration up to $O(q^2)$ with $q$ a parameter controlling the level of coarse graining \cite{KMV, KPRT}. Furthermore, models where only short-range interactions 
appear are inexpensive to simulate, for example with a single spin-flip Metropolis method. 
However, when both short and long-range interactions are present the classical MH algorithm becomes prohibitively expensive
due to the high cost of calculating energy differences arising from the long-range interaction potential. In \cite{KKP} we extend our framework for coupled CGMC to the dynamics case, developing Kinetic Monte Carlo algorithms based on coarse-level rates.
   
%%%%%%%%%%%%%%%%%%%%%%%%%%%%%%%%%%%%%%%%%%%%%%%%%%%%%%%%%%%%%%%%%%%%%%%%%%%%%%%%%%%%%%%%%%%%%
% Sections description
%%%%%%%%%%%%%%%%%%%%%%%%%%%%%%%%%%%%%%%%%%%%%%%%%%%%%%%%%%%%%%%%%%%%%%%%%%%%%%%%%%%%%%%%%%%%%
 
 Section~\ref{theory} describes the classical Metropolis-Hastings   algorithm and some known mathematical theory for convergence and the rate of convergence for MCMC methods.
In Section~\ref{method} we present the proposed Coupled CGMC method in a general framework
  describing  its mathematical properties.   We state the main theorem  that compares the rate of convergence to equilibrium with the rate of the classical MH method.
In Section~\ref{SCG} we describe stochastic lattice systems and the coarse-graining procedure in order to prepare 
for the  application of the proposed method in Section~\ref{example} to a generic  model of   lattice systems 
in which both short and long-range interactions are present.

%
%%%%%%%%%%%%%%%%%%%%%%%%%%%%%%%%%%%%%%%%%%%%%%%%%%%%%%%%%%%%%%%%%%%% 
% 
% MCMC theory and MH
%
%%%%%%%%%%%%%%%%%%%%%%%%%%%%%%%%%%%%%%%%%%%%%%%%%%%%%%%%%%%%%%%%%%%%
%
\section{MCMC methods}
\label{theory}
Before describing the Metropolis-Hastings method we   need to introduce some necessary definitions and theoretical facts.

Let  $X_n$  be a  Markov chain  on space $\Sigma$ with transition kernel $\Kk$.
\begin{defin}
A transition kernel $\Kk$ has the {\it stationary measure} $\mu$ if 
$$\Kk  \mu=\mu\PERIOD$$
 \end{defin}

 \begin{defin}
 $\Kk$ is called reversible with respect to $\mu$ if 

$$ (g,\Kk h)_{\mu}=(\Kk g, h)_{\mu},\ \ \text{ for all } g,h\in L^2(\mu)\PERIOD$$
  \end{defin}
 where $ (g,  h)_{\mu} =\int_{\Sigma} \overline{g(\sigma)}h(\sigma)\mu(d\sigma)$ and $\Kk g(\sigma)= \int_{\Sigma}\Kk(\sigma,d\sigma')g(\sigma'), \forall \sigma \in \Sigma $.
 
 \medskip
 A sufficient condition for $\mu$ being a stationary measure of $\Kk$ is the, often easy to check, detailed balance(DB) condition.
 \begin{defin}
 A Markov chain with transition kernel $\Kk$ satisfies the detailed balance  condition if there exists a function $f$ satisfying 
 \begin{equation}\label{dBD}
 \Kk(\sigma,\sigma')f(\sigma)=\Kk(\sigma',\sigma)f(\sigma') \PERIOD
 \end{equation}
 \end{defin}

Here we focus on the Metropolis-Hastings   algorithm   \cite{RC}.
The algorithm generates an ergodic  Markov chain $X_n$ in the state space $\Sigma$, with stationary measure $\mu(d\sigma)$.
Let $f(\sigma)$ be the probability density corresponding to the measure $\mu$ and $X_0=\sigma_0$ be arbitrary. The $n$-th  
iteration of the algorithm consists of the following steps
\begin{algorithm}[Metropolis-Hastings Algorithm]\label{MH} {\ }  \\

\noindent Given $X_n=\sigma$  
\begin{description}
\item[Step 1] Generate $Y_n=\sigma' \sim q(\sigma'|\sigma)$
\item[Step 2] Accept-Reject
\begin{equation*}
X_{n+1}=
\left\{ {\begin{array}{*{20}l}
   {Y_n=\sigma'\  \text{with probability}\  \alpha(\sigma,\sigma')} \\
   { X_n=\sigma_n\  \text{ with probability} \ 1-\alpha(\sigma,\sigma')} \\
\end{array}} \right.
\end{equation*}
where
\begin{equation*}
\alpha(\sigma,\sigma')=\min\left\{1,\frac{f(\sigma')q(\sigma',\sigma)}{f(\sigma)q(\sigma,\sigma')} \right\} 
\end{equation*}
\end{description}

We denote $q(\sigma'|\sigma) $ the  {\it proposal} probability transition kernel,  
and $\alpha(\sigma,\sigma') $ the probability of accepting the proposed state $\sigma'$. 
The transition kernel associated to MH algorithm is
\begin{equation}\label{KC}
\Kk_c(\sigma,\sigma')=\alpha(\sigma,\sigma')q(\sigma,\sigma')+\left[1-\int\alpha(\sigma,\sigma')q(\sigma,\sigma')d\sigma'\right]\delta(\sigma'-\sigma)\PERIOD
\end{equation}
\end{algorithm}

Convergence and ergodicity properties of the chain $\{X_n\}$ depend on the proposal kernel $q$, and they are studied extensively in \cite{RC}. 
 $\Kk_c$ satisfies the DB condition with $f$ ensuring that it has stationary measure $\mu$.
   $\Kk_c$ is irreducible and aperiodic \cite{RC},
nonnegative definite, and  reversible, thus the Markov chain with transition kernel
$\Kk_c$ converges in distribution to $\mu$.

\subsection{Mixing times and speed of convergence}
\label{mixing}
It is known \cite{DSC} that for a discrete-time Markov chain $X_n$ with the transition kernel 
$\Kk$ and the stationary distribution $f$, the rate of convergence to its stationarity  can be measured in terms of  
kernel's  second largest eigenvalue,  according to 
\[ 
2||\Kk^n(\sigma,\cdot) -f ||_{TV} \le \frac{1}{ f(\sigma)^{1/2}} \beta^n  
\]
where $\beta=max\{|\beta_{min}|,\beta_1 \}$ and $-1\le\beta_{min}\le \dots \le \beta_1\le\beta_0=1 $ are 
the real eigenvalues of $\Kk$. The   spectral gap of kernel $\Kk$ is defined by 
$\lambda(\Kk)=\min \left\{\frac{\mathcal{E}(h,h)}{\VAR(h)}; \VAR(h)\neq 0\right\}$ which for the  self-adjoint kernel $\Kk$,  because of reversibility, is  $\lambda(\Kk) =1-\beta_1$. With the Dirichlet form $\mathcal{E}$ and the variance defined by
\begin{eqnarray*}
&& \mathcal{E}(h,h)=\frac12 \sum_{\sigma,\sigma'}|h(\sigma)-h(\sigma')|^2\Kk(\sigma,\sigma')f(\sigma)\COMMA\\
&& \VAR(h)=\frac12 \sum_{\sigma,\sigma'}|h(\sigma)-h(\sigma')|^2f(\sigma')f(\sigma) \PERIOD
\end{eqnarray*}
Between two algorithms producing   Markov chains with identical equilibrium distributions  {\it better} in terms of the 
speed of convergence is the one with the {\it smaller second eigenvalue in absolute value} or equivalently with the {\it larger  spectral gap}.

\section{The Coupled CGMC method}
\label{method}

The   proposed algorithm is designed to generate samples from the microscopic probability 
measure $\mu$ with density  $f$  on a space  $\Sigma$, coupling properly states of the  microscopic  space $\Sigma$ with states on a  
 {\it  coarse} space  $\bar{\Sigma}$ having less degrees of freedom. A  properly constructed coarse  measure on $\bar{\Sigma}$
will be the basis for constructing efficient  proposal kernels for  MH  algorithms sampling large systems.

The {\it coarsening} procedure is based on  the expansion of 
 the target measure $\mu$ to a coarse and a finer  part.  
 Abstractly we write $f(\sigma)=f(\eta,\xi)$ and $\Sigma=\bar{\Sigma}\times \bar{\Sigma}
'$, where $\eta\in \bar{\Sigma}$ represents the coarse variables.

We denote the projection operator on the coarse variables  
\[ T: \Sigma \to \bar{\Sigma},\ \  T\sigma=\eta \PERIOD\]
The exact coarse marginal is
\[ \bar{f}(\eta)=\int_{\bar{\Sigma}'} f(\eta,\xi)d\xi \PERIOD\]
  To obtain an explicit formula of the coarse marginal is as difficult as sampling the original target distribution since   space $\bar{\Sigma}'$ remains  high dimensional. 
  Therefore use  of  approximating distributions of $\bar{f}$ becomes necessary.
  Such approximations  have been proposed in \cite{KMV1, KPRT} for stochastic lattice  systems and   are abstractly described in Section (\ref{SCG}) and     
for complex macromolecular systems see \cite{briels, kremer, doi, vagelis}.

Denote $\bar{f}_0$ an approximation of  $\bar{f}$ on $\bar{\Sigma}$. 
This distribution, combined with a reconstruction distribution  $f_r(\xi|\eta) $ corresponding to the finer variables $\xi$,   will construct   a candidate for   proposal distribution  in MH algorithms performed  in order to  sample from $f$ at the original space  $\Sigma$. 
An example of a 'good' proposal distribution  is  $f_0(\sigma):=\bar{f}_0(\eta)f_r(\xi|\eta)$. For notational simplicity we write $ f_r(\sigma|\eta)$ instead of $ f_r(\xi|\eta)$.
  In terms of the Metropolis-Hastings algorithm    this means that  $q(\sigma,\sigma')=f_0(\sigma')$, or that $f_0(\sigma')$ is the stationary measure of the proposal kernel $q(\sigma,\sigma')$.
  
 The coupled CGMC  algorithm  is composed of two coupled  Metropolis iterations, the first generating  samples  from the proposal distribution and the second samples from the target measure. 
 The first Metropolis step samples the coarse approximating marginal  $\bar{f}_0(\eta)$, using  an arbitrary proposal transition kernel $\bar{q_0}(\eta,\eta')$ to  produce trial samples $\eta'$. The second step is performed  if the coarse trial sample is accepted,   and consists of  the reconstruction  from the coarse trial state and a   Metropolis criterion designed to ensure  sampling from the correct microscopic density $f$. 
  If a trial coarse sample  is rejected,   then we go back to the first step to rebuild a new coarse trial, so  that the  fine Metropolis step is not performed and no computational time  is wasted on checking fine trial samples that  are most likely to be rejected.

 In \cite{EHL} Efendiev et.al.,  propose  the  Preconditioning MCMC, a two stage ( coarse and fine ) Metropolis MCMC method, applied to  inverse problems of subsurface characterization.  
 The coarse and fine models are finite volume schemes of different resolutions for a PDE two-phase flow model.
 Our algorithm shares the same idea and structure with the Preconditioning MCMC of constructing a proposal density based on meso/macro-scopic properties of the model studied and taking advantage of the first stage  rejections.
   In terms of the MC method 'coarsening'  corresponds to   enriching the range of the sampling measure based on  coarse-scale models proposed by multiscale finite volume methods.
  The major difference of the Preconditioning MCMC and the proposed    algorithm is that  the latter  alternates between
 different state spaces during each MC iteration, the coarse and the finer, whether in the former  the state space remains the same since coarse and fine problems are solved independently. 
 Thus, at the end of a simulation we will have both  fine-scale and "compressed", coarse-grained data.
 The performance of the coarse proposals in our case can be further estimated based on a systematic error analysis such as (\ref{errb}).

The proposed procedure has also some common features with   the modified
  Configurational bias Monte Carlo (CBMS) where the trial density is built up sequentialy with stage-wise rejection decision described in \cite{LIU}, applied effectively  in quantum mechanical systems \cite{CEP}. 
There are also some similarities with simulated sintering and transdimensional MCMC, see  \cite{LIU} and references therein. However, in our method,  the construction of the variable dimensionality (and level of coarse-graining) state spaces and the corresponding Gibbs measures 
relies on statistical mechanics tools that allow a systematic control of the error from one level of coarse-graining to the next, e.g. (\ref{errb}).

\subsection{The algorithm}
\label{algorithm}

\mbox{}
 
 We describe in detail the coupled CGMC Metropolis algorithm outlined in the previews section.
\begin{algorithm}[ Coupled  CGMC   Algorithm]{\ }\\

\noindent 
Let $X_0=\sigma_0$ arbitrary, for $n=0,1,2,\dots$  
 
\noindent 
Given $X_n=\sigma$ 

\begin{description}

\item[Step 1]  Compute the coarse variable $  \eta=T\sigma$

\item[Step 2] Generate a coarse sample $
  \eta'  \sim \bar{q}_0(\eta,\eta')$

\item[Step 3] Coarse Level Accept-Reject

\hspace{22pt} Accept $\eta'$ with probability: 
\[ \alpha_{CG}(\eta,\eta')=\min\left\{1,\frac{\bar{f}_0(\eta')\bar{q}_0(\eta',\eta)}{\bar{f}_0(\eta)\bar{q}_0(\eta,\eta')}\right\}\PERIOD\]

\hspace{22pt}{\bf If $\eta'$ is accepted then proceed to Step 4}

\hspace{22pt}{\bf else generate a new coarse sample Step 2}

\item[Step 4]  Reconstruct $\sigma'$ given  the coarse trial $\eta'$,   
\[\sigma' \sim  f_r(\cdot|\eta')\]

\item[Step 5] Fine Level Accept-Reject 

\hspace{17pt} Accept $\sigma'$ with probability
\begin{align*}
\alpha_f(\sigma,\sigma')=\min\left\{1,\frac{f(\sigma')\bar{f}_0(\eta)f_r(\sigma|\eta)}{f(\sigma)\bar{f}_0(\eta')f_r(\sigma'|\eta')}\right\} \PERIOD
\end{align*}
 
\end{description}
\end{algorithm}

Steps 2 and 3 generate a Markov chain $\{Z_n\}$ in the coarse space $\bar{\Sigma}$  with the transition kernel
 \begin{equation*}\label{Qchain} 
\Qq(\eta,\eta')=\alpha_{CG}(\eta,\eta')\bar{q}_0(\eta,\eta')+\left[1-\int \alpha_{CG}(\eta,z)\bar{q}_0(\eta,z)\right]\delta(\eta'-\eta)\PERIOD
\end{equation*}
The  stationary measure of kernel $\Qq$ is $\bar{f}_0(\eta)$.  
Combination of  this kernel and 
 Steps 1 and 4  constructs the desired  proposal  transition kernel $q_0(\sigma,\sigma')$ on the fine level space $\Sigma$,
\[q_0(\sigma,\sigma')= \Qq(\eta,\eta')f_r(\sigma'|\eta')\PERIOD  \]
According to the MH algorithm in order to sample from $f$, the fine level acceptance probability should be $ \alpha_f(\sigma,\sigma')=\min\left\{1,\frac{f(\sigma')q_0(\sigma',\sigma)}{f(\sigma)q_0(\sigma,\sigma')} \right\} 
  $, but since $\Qq$ satisfies the Detailed Balance   condition $\Qq(\eta,\eta')\bar{f}_0(\eta)=\Qq(\eta',\eta)\bar{f}_0(\eta')$,  $\alpha_f $ is equal to
\begin{align*}
\alpha_f(\sigma,\sigma'
)&=\min\left\{1,\frac{f(\sigma')\Qq(\eta',\eta)f_r(\sigma|\eta)}{f(\sigma)\Qq(\eta,\eta')f_r(\sigma'|\eta')}\right\}\\&=\min\left\{1,\frac{f(\sigma')\bar{f}_0(\eta)f_r(\sigma|\eta)}{f(\sigma)\bar{f}_0(\eta')f_r(\sigma'|\eta')}\right\}\PERIOD
\end{align*}

The chain ${X_n}$ produced by Coupled CGMC algorithm is a Markov chain on the fine space $\Sigma$, with the transition kernel  
\begin{equation}\label{KCG}
\Kk_{CG}(\sigma,\sigma')=\alpha_f(\sigma,\sigma')q_0(\sigma,\sigma')+\left[1-\int\alpha_f(\sigma,\sigma')q_0(\sigma,\sigma')d\sigma'\right]\delta(\sigma'-\sigma)\PERIOD
\end{equation}

The Markov chain  ${X_n}$ generated by the Coupled CGMC algorithm 
converges to the correct  stationary distribution $f$ and is ergodic, which ensures that  $\frac{1}{n} \sum_{j=1}^n h(X_j) $ is a convergent approximation of the averages $ \int h(\sigma) f(\sigma)d\sigma$ for any $h\in L^1(f)$.
Ergodicity and reversibility properties are satisfied ensuring that the algorithm generates samples from the correct measure.We state this fact as a separate theorem proof of which is given in detail in \cite{KKP}.

We denote 
$E=\{\sigma\in \Sigma; f(\sigma)>0 \} $, $\bar{E}=\{\eta\in \bar{\Sigma}; {\bar f}_0(\eta)>0 \} $.
 \begin{thm}
 For every conditional distribution $\bar{q}_0 $, and $f_r$ such that the   support of $q_0f_r$ includes $E$,
 
 \begin{enumerate}[i)]
 \item  The transition kernel satisfies the detailed balance (DB) condition with $ f $
$$\Kk_{CG}(\sigma,\sigma')f(\sigma)=\Kk_{CG}(\sigma',\sigma)f(\sigma') $$
 \item    $f$ is a stationary distribution of the chain.
 \item if $q_0(\sigma,\sigma')>0,\ \forall \sigma,\sigma'\in E$ and $E \subseteq supp(f_0) $ then $X_n$ is  $f$-irreducible   
 \item is aperiodic  
 \end{enumerate}
 \end{thm}

\subsection{The rate of convergence}
\label{rate}
 
  The calculation of the rate of convergence to stationarity is a hard problem   since it is model dependent as argued earlier. 
  What we can prove though for the proposed method is that it is comparable to the classical Metropolis-Hastings algorithm described in Algorithm \ref{MH}. This fact is stated rigorously in the following theorem which we prove in \cite{KKP}.

 Let  $\lambda(\Kk_{CG}), \lambda(\Kk_{c})$ be the spectral gap corresponding to  the coupled  CGMC $\Kk_{CG}$, (\ref{KCG}), and the classical MH  
  $\Kk_c $, (\ref{KC}),  transition kernels respectively.
 
 \begin{thm}\label{Theo}
Let $q(\sigma,\sigma') $ be  a symmetric proposal transition probability for the classical MH algorithm and $\bar{q}_0(\eta,\eta') $  a symmetric proposal transition probability  on the coarse space $ \bar{\Sigma}$ for the coupled CGMC algorithm, then for any reconstruction conditional  probability $f_r(\sigma|\eta)$

\begin{enumerate}[i)]
\item \begin{equation}\label{kernel}
  \Kk_{CG}(\sigma,\sigma')=\mathcal{ A}(\sigma,\sigma')\mathcal{B}(\sigma,\sigma') \Kk_{c}(\sigma,\sigma')\end{equation}

\medskip
 \begin{equation*}
\mathcal{B}(\sigma,\sigma') =\left\{  \begin{array}{cc} 
\frac{\bar{q}_0(\eta,\eta') f_r(\sigma'|\eta')}{q(\sigma,\sigma')}\\
\frac{\bar{q}_0(\eta',\eta) f_r(\sigma|\eta)}{q(\sigma',\sigma)}
\end{array}  \right.
\end{equation*}

Furthermore we define the subsets
\begin{align*}
&C_1=\left\lbrace (\sigma,\sigma')\in \Sigma\times\Sigma: \left\lbrace \alpha<1, \alpha_{CG}<1, \alpha_f<1 \right\rbrace \  \text{ or }\   \left\lbrace \alpha=1, \alpha_{CG}=1, \alpha_f=1  \right\rbrace \right\rbrace\\
&C_2=\left\lbrace  (\sigma,\sigma')\in \Sigma\times\Sigma:\left\lbrace \alpha=1, \alpha_{CG}<1, \alpha_f=1 \right\rbrace\  \text{ or }\   \left\lbrace \alpha<1, \alpha_{CG}=1, \alpha_f<1  \right\rbrace \right\rbrace\\
&C_3=\left\lbrace  (\sigma,\sigma')\in \Sigma\times\Sigma: \left\lbrace \alpha=1, \alpha_{CG}=1, \alpha_f<1  \right\rbrace\  \text{ or }\   \left\lbrace  \alpha<1, \alpha_{CG}<1, \alpha_f=1  \right\rbrace \right\rbrace\\
&C_4=\left\lbrace  (\sigma,\sigma')\in \Sigma\times\Sigma:  \left\lbrace\alpha<1, \alpha_{CG}=1, \alpha_f=1 \right\rbrace\  \text{ or }\   \left\lbrace  \alpha=1, \alpha_{CG}<1, \alpha_f<1 \right\rbrace \right\rbrace
\end{align*}

\begin{equation*}
\mathcal{A}(\sigma,\sigma')=
\left\{  \begin{array}{cc}
    1,\ \  &\text{ if } (\sigma,\sigma')\in  C_1 \\
   \min\{ \frac{\bar{f}_0(\eta')}{\bar{f}_0(\eta)},\frac{\bar{f}_0(\eta)}{\bar{f}_0(\eta')}\} , \ \  &\text{ if } (\sigma,\sigma')\in  C_2  \\
     \min\{ \frac{f(\sigma')\bar{f}_0(\eta)}{f(\sigma)\bar{f}_0(\eta')},  \frac{f(\sigma)\bar{f}_0(\eta')}{f(\sigma')\bar{f}_0(\eta)}\} , \ \  &\text{ if } (\sigma,\sigma')\in  C_3
      \\
     \min\{ \frac{f(\sigma')}{f(\sigma)},  \frac{f(\sigma)}{f(\sigma')}\} , \ \  &\text{ if } (\sigma,\sigma')\in  C_4
\end{array}  \right.
\end{equation*}

 \item  
 \begin{equation}\label{gap}  \mathcal{A}\underline{\gamma}\lambda(\Kk_c) \le \lambda(\Kk_{CG}) \le \bar{\gamma}\lambda(\Kk_c)
\end{equation}

\medskip
where $\mathcal{A}=\inf_{\sigma,\sigma'}\mathcal{ A}(\sigma,\sigma')$
and $\underline{\gamma}>0, \bar{\gamma}>0 $ such that 
$  \underline{\gamma}\le \Bb(\sigma,\sigma')\le \bar{\gamma}\PERIOD$
 \end{enumerate}
 \end{thm}

%%%%%%%%%%%%%%%%%%%%%%%%%%%%%%%%%%%%%%%%%%%%%%%%%%%%%%%%%%%%%%%%%%%%%%  
% Extended Lattice Systems, Coarse graining and Reconstruction
%%%%%%%%%%%%%%%%%%%%%%%%%%%%%%%%%%%%%%%%%%%%%%%%%%%%%%%%%%%%%%%%%%%%%%

\section{Extended Lattice Systems}\label{SCG}

 This class of stochastic processes   is employed in the modeling of
 adsorption, desorption, reaction and
diffusion of chemical species 
in numerous applied science areas such as catalysis, microporous materials,  biological systems, etc. \cite{au, binder}. 
To demonstrate the basic ideas, we
consider an Ising-type system
on a periodic $d$-dimensional lattice $\LATT$ with $N=n^d$ lattice points.
At each 
$x \in \LATT$ we can define an 
order parameter $\sigma(x)$; for instance, when taking    values 
$0$ and $1$, it can   describe vacant and occupied sites.
The energy $H_N$ of the system,
 at the configuration $\sigma=\{ \sigma (x):x\in\LATT\}$ is given by the Hamiltonian,
\begin{equation}
H_N(\sigma)=-{1 \over 2}\sum_{x \in \LATT}\sum_{y\not= x}\left[K(x-y) + J(x-y)\right]\sigma (x)\sigma
(y)+\sum h \sigma (x)\COMMA 
\label{hamiltonian}
\end{equation}
where $h$, is the  external field and $J$
is the inter-particle  potential.
Equilibrium states  at the temperature $\sim \beta^{-1}$ are
described by the (canonical) Gibbs probability  measure, and  $Z_{\LATT}$ is the normalizing constant (partition function)
\begin{equation}\label{gibbs}
\mu_{\LATT, \beta}(d \sigma)
= Z_{\LATT}^{-1}\exp\big(-\beta H_N(\sigma)\big)P_N(d\sigma)\PERIOD
\end{equation}
Furthermore, the product Bernoulli distribution  $P_N(\sigma)$  is  the  {\it prior distribution} on $\LATT$.
 
 The inter-particle potentials $ K,\ J$  account for  interactions between occupied sites.
We consider $K$ corresponding to the short and $J$ to the long-range interactions discussed in detail in Section (\ref{lattice-sl}). General potentials with combined short and long-range interactions are discussed here, while we can also  address  potentials   with  suitable decay/growth  conditions \cite{AKPR}.

The prior  $P_N(d\sigma)$ is typically a product measure, describing  the system at $\beta=0$, when  interactions in $H_N$ are unimportant
and thermal fluctuations-disorder-associated with the product structure of $P_N(d\sigma)$ dominates. By contrast at zero temperature, $\beta=\infty$ interactions, and hence order, prevail. Finite temperatures, $0<\beta <\infty$,  describe intermediate states, including  possible phase transitions between ordered and  disordered states.
 For both on-lattice or off-lattice particle  systems, the finite-volume equilibrium states of the system have the structure   \eqref{gibbs}.

\subsection{ Coarse-graining of microscopic   systems}

Coarse-graining of microscopic   systems is essentially an approximation theory and a numerical analysis question.  However,  the presence of {\em stochastic fluctuations} on one hand, and the {\em extensivity} of the models (the system size scales
with the number of particles) on the other, create a new set of challenges. We discuss all these issues next, in a general setting that applies to both on-lattice and off-lattice systems. 

First,  we write the microscopic configuration $\sigma$ in terms of coarse variables $\eta$ and  corresponding  fine ones  $\xi$ so that $\sigma=(\eta, \xi)$. 
We denote  
by $T$ the coarse-graining map 
$T\sigma=\eta$.
 
The CG system size is denoted by $M$, while the microscopic system size is $N=Mq$,
where we  refer to $q$ as the level of coarse graining,
  and $q=1$ corresponds to no coarse graining.
The  exact   CG Gibbs measure  is given (with a slight abuse of notation) by
$
\BARIT{\mu}_{M,\beta} \,=\,  \mu_{N,\beta} \circ  T^{-1} \,.
$
In order to write $\BARIT{\mu}_{M,\beta}$ in a more convenient form
we first define the CG prior
$ 
\BARIT{P}_M( d{\eta}) = P_N \circ T^{-1} 
%= \prod_{k \in \LATTC} \bar{P}_k (d \eta(k))
%\COMMA \quad {\rm with} \quad \BARIT{\rho}(\eta(k) ) \,=\,
%\binom{q}{\frac{\eta(k)+ q}{2}}\left(\frac{1}{2}\right)^{q} \PERIOD
$.   
The  conditional  prior  $P_N(d\sigma|{\eta})$ is the  probability of having 
a microscopic configuration $\sigma$,  given a coarse configuration $\eta$.
We now rewrite  $\BARIT{\mu}_{M,\beta}$ using the {\em exact coarse-grained 
Hamiltonian}:
\begin{eqnarray}\label{rg}
     \EXP{-\beta \BARH_M (\eta)} & \,= \,&\EXPECT[\EXP{-\beta H_N}|{\eta}]=\int \EXP{-\beta H_N(\sigma)} P_N(d\sigma|{\eta})\COMMA
\end{eqnarray}
a procedure  known as the {\em renormalization group map},  \cite{Golden}; 
  $\BARIT{\mu}_{M,\beta}(d{\eta})$  is now re-written using \eqref{rg} as 
\begin{equation}\label{cg_gibbs}
     \BARIT{\mu}_{M,\beta}(d{\eta})=\frac{1}{\BARIT{Z}_M}
     \EXP{-\beta\BARH_M({\eta})} \BARIT{P}_M( d{\eta})\PERIOD
\end{equation}

Although typically $\BARIT{P}_M( d{\eta})$ is  easy to calculate,
even for moderately small values of $N$ the exact computation of the coarse-grained 
Hamiltonian $\BARH_M(\eta)$ given by \eqref{rg} is, in general, impossible.

We have  shown in \cite{KPRT} that there is an expansion of   $\BARH_M(\eta)$ into a convergent series 
\begin{equation}
\label{series}
\BARH_M(\eta)=\BARH_M^{(0)}(\eta)+\BARH_M^{(1)}(\eta)+\BARH_M^{(2)}(\eta)+\cdots+\BARH_M^{(p)}(\eta) + N \times \BIGO(\epsilon^{p})\, ,
\end{equation}
by constructing  a suitable  first approximation $\BARH_M^{(0)}(\eta)$ and identifying a suitable small 
parameter $\epsilon$  to control the higher order terms in the expansions. Truncations including the first terms in  \eqref{series}
correspond to coarse-graining schemes of increasing accuracy.
In order to obtain this expansion  we  
rewrite  \eqref{rg} as
\begin{equation}\label{rg2}
\BARH_M(\eta)=\BARH_M^{(0)}(\eta)
-\frac{1}{\beta}\log \EXPECT[\EXP{-\beta (H_N-\BARH_M^{(0)}(\eta))}|{\eta}]\PERIOD
\end{equation}
We need    to show that the logarithm can be expanded into a convergent series, yielding eventually  \eqref{series}, however, two interrelated  difficulties emerge immediately: first,   the stochasticity of the system in the finite temperature case, yields  the nonlinear log expression  
which in turn will need to be expanded into a series. Second,  the extensivity of the microscopic system, i.e., typically the Hamiltonian scales as $H_N=O(N)$, does not allow the expansion of the logarithm and exponential functions into a Taylor series.
For these two reasons, one of the  mathematical tools we employed  
is the   {\em cluster expansion method}, see  \cite{Simon}  for an overview. Cluster expansions   allow us to identify uncorrelated components in the expected value
$\EXPECT[\EXP{-\beta (H_N-\BARH_M^{(0)}(\eta))}|{\eta}]\COMMA$
which in turn will permit us to  factorize it,
and subsequently expand the logarithm. 

%%%%%%%%%%%%%%%%%%%%%%%%
\medskip
The coarse-graining of systems with purely long- or intermediate-range
interactions of the form
\begin{equation}\label{microJL}
J(x-y)=L^{-1}V\Big({(x-y)/ L}\Big)\, , \quad x\, , y \in \LATT\, ,
\end{equation}
where
$V(r)=V(-r)$, 
$
V(r)=0\, , |r| > 1
$, 
was studied using cluster expansions  in \cite{KPRT, AKPR, KPR}. 
The corresponding   
  CG Hamiltonian is 
\begin{equation}\label{Hzero}
\bar H^0(\eta)=
-{1 \over 2}\sum_{l \in \LATTC}\sum_{k \in \LATTC
\atop k\ne l}\bar J(k, l)\eta(k)\eta(l)
-
{\bar J(0, 0) \over 2}\sum_{l \in \LATTC}\eta(l)\big(\eta(l)-1\big)
+\sum_{k \in \LATTC}\bar h
\eta(k)\, .
\end{equation}

\bigskip
$$\bar{J}(k,l)=\frac{1}{q^2}\sum_{x\in C_k}\sum_{y\in C_l} J(x-y),\ \ \bar{J}(k,k)=\frac{1}{q(q-1)}\sum_{x\in C_k} \sum_{y\in C_k,y\neq x}J(x-y)$$ 

%%%%%%%%%%%%%%%%%%%%%%%%
  
One of the results therein is on deriving  error estimates  in terms of the specific relative
entropy $\Rr(\mu| \nu):=N^{-1}\sum_\sigma \log\big\{\mu(\sigma)/
\nu(\sigma) \big\}
\mu(\sigma) \quad 
$ between  the corresponding equilibrium Gibbs measures. Note that the scaling factor $N^{-1}$ is related to the extensivity of the system, hence the proper error quantity that needs to be tracked is the 
loss of information {\it per particle}. Using this idea we can 
assess the {\em  information compression}  for the same level of coarse
graining in  schemes differentiated  by
the truncation level $p$ in (\ref{series})
\begin{equation} 
 \Rr(\BARIT{\mu}_{M,\beta}^{(p)}|\mu_{N,\beta}\circ T^{-1}) = 
    \BIGO(\epsilon^{p+1})\COMMA \quad\quad \epsilon \,\equiv\,    \beta  \|\nabla V\|_1 \big(\frac{q}{L}\big)\, ,\label{errb}
\end{equation}
where $\BARH_M^{(0)}(\eta)$ in \eqref{series} is given by \eqref{Hzero}.
The role of such higher order schemes was demonstrated 
in nucleation, metastability and the resulting  switching times between phases, \cite{AKPR}. 
 
Although CGMC and other CG methods can provide a powerful computational tool in molecular simulations,
it has been observed  that in some regimes,  important macroscopic properties may not be  captured properly.
For instance, (over-)coarse graining in  polymer systems may yield  wrong predictions   in the melt structure 
\cite{AKr};
similarly wrong predictions on    crystallization were also observed in the 
CG of complex fluids, \cite{PK06}. 
In CGMC  for lattice systems,
hysteresis and critical behavior may also  not be captured properly
for short and intermediate range potentials, \cite{KMV, KPRT}. 
Motivated by such observations,
in our recent work we studied when CG  methods perform satisfactorily, and how to quantify the CG  approximations  from a {\em numerical analysis} perspective, where error is assessed in view of a specified tolerance. 
Next, we discuss systems with {\em 
long range} interactions, i.e.,    $L>>1$  in (\ref{microJL}). These systems  can exhibit  complex behavior such as phase transitions, nucleation, etc.,  however, they are  more tractable analytically. At the same time   they pose a serious challenge to conventional MC methods due to the large number of neighbors involved in each MC step.

Here we adopt this  general approach, however, the challenges when both short and long-range interactions are present,  require a new methodology. Short-range interactions induce strong "sub-coarse grid"    fine-scale correlations between coarse cells, and need to be explicitly included in the initial approximation $\BARH_M^{(0)}(\eta)$. For this reason we introduced in \cite{KPRT3} a {\em multi-scale decomposition}  of the Gibbs state \eqref{gibbs},
 into fine and coarse variables, which in turn  allows us to describe in an explicit manner the communication across scales, for both short and long-range interactions.

\subsection{Multiscale Decomposition and Splitting Methods for MCMC}\label{lattice-sl}

We first focus on  general lattice systems, and subsequently discuss related applications in later sections. 
We
 consider   (\ref{hamiltonian}) where in addition to the long-range potential  \eqref{microJL}, we add the short-range 
$K(x-y) = S^{-1}U\left(N|x-y|/S\right)
$,
where $S<<L$ and $U$ has similar properties as $V$ in (\ref{microJL}); 
for $S=1$ we have the  usual nearest neighbor interaction. 
The  new Hamiltonian includes both long and  short-range interactions: 
$H_N=H_N^{l}+H_N^{s}\, .$  

{\em The common theme}  is the observation that long-range interactions $L>>1$  can be handled very efficiently by CGMC, \eqref{errb}.
On the other hand short-range interactions are relatively inexpensive and one could  simulate them with Direct Numerical Simulation (DNS)  provided there is a suitable {\em splitting}
of the algorithm in short and long-range parts, that  can reproduce within a given tolerance equilibrium  Gibbs states and  dynamics.
We return to the general discussion in \eqref{series} and outline the steps we need  in order to construct the CG Hamiltonian for the combined short and long-range interactions.

\medskip
\noindent
{\em Step 1: Semi-analytical splitting schemes.} Here we take advantage of CG approximations developed in \eqref{errb} in order to decompose our calculation into  analytical and   numerical components, the latter involving   only 
short-range interactions:
\begin{eqnarray*}
\mu_{N,\beta}(d\sigma) &\sim& e^{-\beta H_N^{} (\sigma)  } P_N(d\sigma) =\\&=& e^{-\big(\beta H_N^{l}(\sigma)-\bar H_M^{l,0}(\eta)\big)}\Big[e^{-\beta H_N^{s} (\sigma)  } P_N(d\sigma | \eta)\Big]
e^{-\bar H_M^{l,0}(\eta)}\bar P_{M}(\eta)\COMMA
\end{eqnarray*}
where $\bar H_M^{l,0}$ is the analytical CG formula \eqref{Hzero} constructed for the computationally expensive, for conventional MC, long-range part; due to the estimates \eqref{errb}, the first term 
has controlled error. Furthermore, the dependence of $\epsilon$ on $\nabla V$ in these  estimates    suggests  a {\em rearrangement} of  the overall  combined short- and long-range potential, 
into  a new short-range interaction   that includes  possible singularities originally  in the long-range component \eqref{microJL}, e.g., the singular  part in a Lennard-Jones potential,
and  a locally integrable (or smooth)  long-range decaying component that can be analytically coarse-grained using \eqref{Hzero}, with a small error due to \eqref{errb}. This breakdown allows us to isolate the short-range interactions (after a possible re-arrangement!), 
and   suggests the two alternative computational approaches:  either seek an approximation $e^{-\beta   \BARH_M^{s}(\eta)}=\int e^{-\beta H_N^s} P_N(d\sigma|\eta)$, or  use sampling methods  
to account for the short-range  "unresolved" terms.   

\subsection{Microscopic Reconstruction}\label{recon}
The reverse procedure of coarse-graining, i.e., reproducing "atomistic" properties,  directly from CG simulations is an issue that arises extensively  in the polymer science literature, \cite{tsop, mp}. The principal idea  is that computationally inexpensive CG simulations 
will reproduce the large scale structure and subsequently microscopic information will  be added through 
{\em microscopic reconstruction}, e.g., 
the calculation of  diffusion of penetrants through polymer melts, reconstructed from  CG simulation, \cite{mp}.
In this direction,   CGMC provides a simpler lattice framework to mathematically formulate microscopic  reconstruction  and study related numerical  and computational issues. Interestingly this issue  arised also in the mathematical error analysis in \cite{KT, KPS}.

The mathematical formulation for the reconstruction of the microscopic 
equilibrium  follows trivially when   we rewrite  the Gibbs measure \eqref{gibbs} 
in terms of   the exact CG measure corresponding to \eqref{rg}, defined in \eqref{cg_gibbs}, \cite{KPR}:
$${\mu}_{N}(d \sigma) \sim e^{- \beta( H(\sigma) - \bar H(\eta))} P_N( d\sigma \vert \eta) 
\bar{\mu}_{M}(d \eta)  \equiv\, \mu_N(d\sigma\vert \eta) \bar{\mu}_{M}(d \eta) \,.
$$
We can define   the conditional probability $\mu_N(d\sigma\vert\eta)$ as the {\em exact reconstruction} of 
$\mu_N(d\sigma)$ from  the exactly CG measure $\bar \mu_M(d \eta)$. Although many fine-scale configurations 
$\sigma$ correspond to a single CG configuration $\eta$,
the ``reconstructed" conditional probability $\mu_N(d\sigma\vert\eta)$ is {\em uniquely} defined, 
given the microscopic and the coarse-grained measures $\mu_N(d\sigma)$ and  $\bar \mu_M(d \eta)$ respectively.

A coarse-graining scheme provides  an approximation $\bar \mu_{M}^{\rm app}(d \eta)$ 
for $\bar \mu_M(d \eta)$,  at the coarse level.  The approximation  $\bar \mu_{M}^{\rm app}(d\eta)$ could be, for instance,   any of the schemes discussed in Section  \ref{lattice-sl}. 
To provide a reconstruction we need to lift the measure $\bar \mu_{M}^{\rm app}(d \eta)$ to
a measure $\mu_{N}^{\rm app}(d \sigma)$ on the microscopic configurations. That is,  we  
need to specify a conditional probability $\nu_N(d\sigma \vert  \eta)$ and set
$\mu_N^{\rm app}(d\sigma)\,:=\, \nu_N(d\sigma \vert  \eta) \bar \mu_{M}^{\rm app}(d\eta)\,.
$
In the spirit of our earlier discussion, it is natural to measure the efficiency of the reconstruction  by the  relative entropy, 
\begin{equation}\label{entdec}
{\cal R} \left( \mu_N^{\rm app}  \SEP {\mu}_{N} \right)\,=\, {\cal R} \left( \bar \mu_M^{\rm app}  \SEP \bar {\mu}_{M} \right) + 
\int  {\cal R} \left( \nu_N( \cdot \vert \eta) \SEP \mu_N( \cdot \SEP \eta) \right)  \bar{\mu}_{M}^{\rm app}(d \eta) 
\,,
\end{equation}     
i.e., relative entropy splits the total error at the microscopic level into the sum of the error at the coarse level and the 
error made during  reconstruction, \cite{KPR, TT}.  

The first term in \eqref{entdec} can be controlled via CG estimates, 
e.g.,  \eqref{errb}.   However, \eqref{entdec} suggests that in order to obtain a successful reconstruction we then need to construct $\nu_N(d\sigma\,\vert\,\eta)$ such that (a) ${\cal R} \left( \nu_N( d\sigma\,\vert\, \eta) \SEP \mu_N( d\sigma\SEP \eta) \right)$ 
should be of the same order as the first term in  \eqref{entdec}, and (b) it is easily computable and implementable.

\medskip
\noindent
 
The simplest example of reconstruction is obtained by considering a microscopic system with intermediate/long-range interactions (\ref{microJL})
\begin{equation}\label{recon0a}
\bar{\mu}_M^{\rm app}(d \eta) \,=\, \bar{\mu}_{M}^{(0)} ( d\eta) \,, \quad  \nu_N(d\sigma\,\vert\, \eta) = 
P_N(d\sigma\,\vert\,\eta) \PERIOD
\end{equation}
Thus we first sample the CG variables $\eta$ involved in $\bar{\mu}_{M}^{(0)}$, using a CGMC algorithm; then  we reconstruct the microscopic configuration $\sigma$ by distributing the particles uniformly on the 
coarse cell, conditioned on the value of $\eta$. Since $P_N(d\sigma \vert \eta)$ is a product measure
this can be done numerically in a very easy way, without communication between coarse cells and only at the coarse cells where an update has occurred in the CGMC algorithm.  In this case the analysis in \cite{KPRT2} yields the estimates
$$
{\cal R} \left( \bar \mu_M^{(0)} \, \vert \, \bar {\mu}_{M} \right) = O(\epsilon^2) \,,\,\,\,   
{\cal R} \left( \mu_N( \cdot \SEP \eta) \,\vert\, P_N( \cdot \SEP \eta) \right) \,=\,  
\frac{\beta}{N}\left( \bar H^{(0)}(\eta)- \bar H(\eta) \right)\,=\, O(\epsilon^2)\PERIOD
$$
Hence the reconstruction is second order accurate and  of the same order as the coarse-graining given by \eqref{Hzero}.

%%%%%%%%%%%%%%%%%%%%%%%%%%%%%%%%%%%%%%%%%%%%%%%%%%%%%%%%%%%%%%%%%%%%  
%Examples 
%%%%%%%%%%%%%%%%%%%%%%%%%%%%%%%%%%%%%%%%%%%%%%%%%%%%%%%%%%%%%%%%%%%%
%
\section{Example: Short and long-range interactions}\label{example}
Short and long-range interactions pose a formidable computational challenge. 
We consider  an example that has been  explicitly solved by Kardar in \cite{KAR}.
The model considered has state space $\Sigma_N = \{0,1\}^{\Lambda_N}$, where $\Lambda_N$ is a 1-dimensional lattice  with $N$ sites.
The energy of the system at   configuration $\sigma=\{\sigma(x), x\in\Lambda_N\} $ is

\begin{eqnarray*}
\beta H_N(\sigma) &=& -\frac{K}{2}\sum_{x}\sum_{|x-y|=1}\sigma(x)\sigma(y) -\frac{J}{2N}\sum_{x}\sum_{y\neq x} \sigma(x)\sigma(y)-h \sum \sigma(x)\\
&\equiv&  H^s_N(\sigma)+ H^l_N(\sigma)+E(\sigma)\PERIOD
\end{eqnarray*}
  Hamiltonian $H_N(\sigma)$ consists of the short-range  term $H^s_N$,  the long-range term $H^l_N $ and 
an external field $E$.
The interactions involved in $H^s_N$ are of  the nearest-neighbor type with strength $K$, while $H^l_N $  represents
a  mean-field approximation or the Curie-Weiss model defined by the potential $J$ averaged over all lattice sites.  
For this generic model Kardar gave in \cite{KAR} a closed form solution for magnetization $M_{\beta}(K,J,h)$,for  the state space $\{-1,1\}$
 
$$M_{\beta}(K,J,h)=\arg\min_{m}\left( \frac{J}{2}m^2-\log\left[ e^K\cosh(h+Jm)+\sqrt{e^{2K}\sin^2(h+Jm)+e^{-2K} }\right]\right) $$  
a simple rescaling of which gives the exact average coverage $m_{\beta}(K,J,h)$ for the lattice-gas model considered here.
 
 \begin{equation}\label{cov}
 m_{\beta}(K,J,h) =\frac12\left( M_{\beta}\left(\frac14 K,\frac14 J,\frac12 h - \frac14 J - \frac14 K\right) +1\right)
 \end{equation}
 
 We have constructed the classical single spin-flip M-H algorithm and the coupled Metropolis CGMC for the single spin-flip algorithm, both generating samples from the  Gibbs measure $$\mu_{N,\beta}=\frac{1}{Z_N}e^{-\beta H_N(\sigma)}P_N(d\sigma)\PERIOD $$
 
We denote $\sigma^x$ the state that differs from $\sigma$ only at the site $x$,   $\sigma^x(y)=\sigma(y),y\neq x $, $\sigma^x(x)=1-\sigma(x)$, the  proposal transition kernel is $q(\sigma'|\sigma) = \frac{1}{N}\sum_{x}\delta(\sigma'-\sigma^x)  $, proposing a spin-flip at the site $x$ with the probability $\frac1N$. 

We apply the coupled   CGMC  method  with coarse updating  variable 
$$\eta:=T\sigma=\{\eta(k), k=1,\dots,  M\}$$\\ 
$\eta(k):=\sum_{x\in C_k} \sigma(x),\ qM=N $ with a coarsening level $q<M$,
 and for  the  maximum coarsening    $q=N$ where the coarse variable is  total magnetization $\eta=\sum_{x\in \Lambda_n} \sigma(x)$.
This can be thought as a coarsening procedure
constructing  a system consisting of one big coarse cell $M=1$ with $q=N$ sites.
  Since we want to consider only single spin-flip updates, for the sake of comparison to the classical Metropolis method,  the cell updating 
can take only the values $\pm 1$ and the reconstruction is chosen uniform in each cell, in the  sence described in example at Section \ref{recon}.
 
 Table \ref{cost2}  gives a comparison of  
the classical single-site updating  Metropolis Hastings algorithm   with the proposed coupled Metropolis CGMC algorithm,
 in terms of  computational complexity per iteration. By computational complexity here  we mean the cost of calculating energy differences involved at the acceptance probabilities. 
 Consider the case that both the microscopic single-site updating Metropolis and the two-step CGMC are run $n$ times. This is reasonable to consider since   as  stated   at Theorem \ref{Theo}   the two methods have comparable mixing times, therefore the number of iterations needed to achieve stationarity are  comparable.
 We denote $\E(\alpha_{CG}):=\int\int \alpha_{CG}(\eta,\eta')\bar{q}_0(\eta,\eta') \bar{f}_0(\eta)d\eta d\eta' $    the average acceptance rate of the coarse proposal. The average number of accepted coarse samples
  is $ n_1:=[\E(\alpha_{CG}) n]$ for which $n_1<n $ since $\E(\alpha_{CG})< 1$ .
 This means  that the reconstruction and fine step acceptance criterion are performed  in average only for $n_1$ iterations.

\begin{table} 
 \caption{Operations count for  evaluating energy differences  for $n$ iterations }
  \begin{tabular}{lccc} 
\hline\noalign{\smallskip}
 Cost & Metropolis Hastings  & Coupled CGMC $q<N$ & Coupled CGMC $q=N$\\
\noalign{\smallskip}\svhline\noalign{\smallskip}
coarse A-R & -- &  $n \times  O(M)$  & $n \times  O(1)$\\
fine A-R &  $n \times O(N)$ & $n_1 \times O(1)$ & $n_1 \times O(1)$ \\
\noalign{\smallskip}\hline\noalign{\smallskip}
 \end{tabular}
 \label{cost2}
 \centering
 \end{table}

\begin{figure}[h]
	 \centering
	 \includegraphics[scale=0.4]{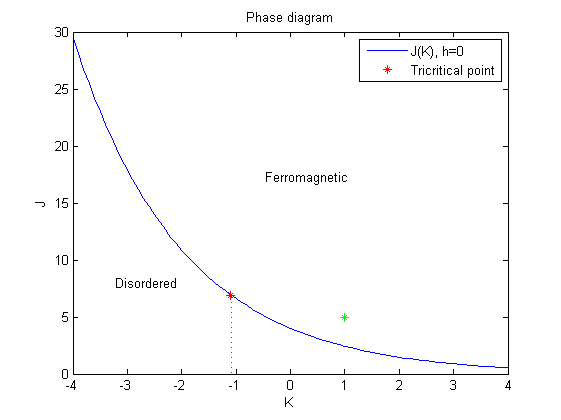}
	    
	 \caption{  Phase Diagram \cite{KAR}  }
\label{fig:1}        
\end{figure}

\begin{figure}[h]
	\subfigure{\includegraphics[angle=-0,width=0.45\textwidth, height=0.3\textheight]{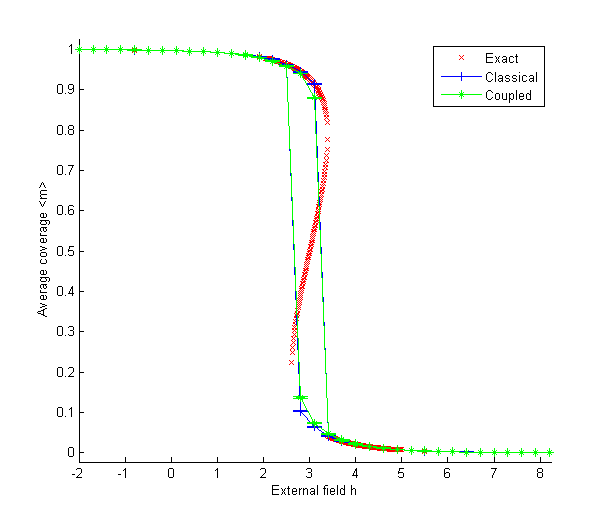}}
	        \hspace{.5cm}
	\subfigure{\includegraphics[angle=-0,width=0.45\textwidth, height=0.3\textheight]{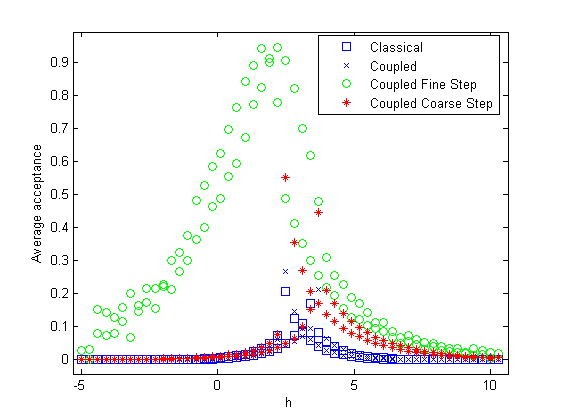}}
	 \caption{N=  1028, q=8, 
K= 1,
J= 5: (a) Coverage  ;  (b) Average acceptance.}
\label{fig:2}
\end{figure}

 \begin{figure}[h]
	 \centering
	 \includegraphics[scale=0.4]{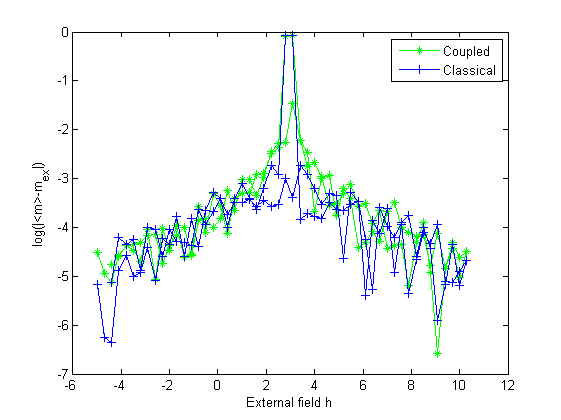}
	    
	 \caption{N=  1028, q=8, 
K= 1,
J= 5: Local error $ \log(|<m>-<m_{ex}|)$   }
\label{fig:3}        
\end{figure}

Results of computational implementation are shown in Figure \ref{fig:2} and Table \ref{cost}  and \ref{cost1}. Figure \ref{fig:2}a represents the average coverage  versus the external field $h$ for the exact solution $m_{ex} $, the classical MH result $ <m_{cl}>$  and the coupled CGMC  $<m>$, for a choice of interaction parameters $K=1,\ J=5$ in the ferromagnetic region as is stated at the phase diagram depicted in Figure \ref{fig:1}. The exact solution $ m_{ex}$ as is ploted in Figure  \ref{fig:2}a corresponds to the part the full solution (\ref{cov}) up to the point it jumps.
 Figure \ref{fig:2}b is a graph of the average acceptance rates for the classical MH algorithm
   and the coupled CGMC algorithm, that verifies the theoretical proof of 
 the fact that the two algorithms have comparable mixing times since the acceptance rate is strongly related to mixing times. In the same figure we also give the average acceptance rates of the coarse   and fine step   of the coupled method, noting that the fine acceptance rate is high which means that  most 
 of the trial samples entering the fine step are accepted.

\begin{table}[ht]
\caption{N= 4096}  
\centering  
\begin{tabular}{l c c c }  
\hline\noalign{\smallskip} 
  & CG level q &  $\text{Error}_{c}$ & CPU(min)
\\ [0.5ex]
 \noalign{\smallskip}\svhline\noalign{\smallskip}
& 4& 0.089 &   93.52       \\[-1ex]
\raisebox{1.5ex}{ K= -2.0, J=2 } &
  8 & 0.302     & 45.8 \\[1ex]
  &4 & 0.003   & 93.6   \\[-1ex]
\raisebox{1.5ex}{K= 1.0, J=5 } & 
8 &  0.003  &  45.9 \\[1ex]
  &4 & 0.027  & 91.6   \\[-1ex]
\raisebox{1.5ex}{K= 1, J=1} & 
8&  0.100 &    45.5 \\[1ex]
\hline  
\end{tabular}
\label{cost}
\end{table}

Table \ref{cost}
reports the error between the exact solution and the average coverage obtained from the coupled CGMC algorithm. Error is measured in terms of the pointwise solutions as $\text{Error}_{c}=\left(\sum_{i}(m_{ex}(h_i)-<m>(h_i))^2\right)^{1/2} $ and $\text{Error}_{cl}=\left(\sum_{i}(m_{ex}(h_i)-<m_{cl}>(h_i))^2\right)^{1/2} $ for the coupled and the classical method respectively, where $h_i$ are the different external field parameters for which the average coverages are computed.
CPU times are compared for the coarse-graining levels $q=4$ and $q=8$. To demonstrate the robustness of the algorithm we present simulations at three different points of the phase diagram plane $K-J$: in the disordered (($K= -2.0, J=2$) and ($K= 1, J=1)$) and   ferromagnetic ($K= 1.0, J=5$) regions. 
In table \ref{cost1} we compare  the error between the coupled CGMC average 
coverage with the exact solution and the corresponding CPU time for $q=4$ and $q=8$, in the ferromagnetic region ($K= 1.0, J=5$) for which the classical Metropolis results.

These results demonstrate the efficiency of the  coupled CGMC methods in terms of computational time since the run time gain scales almost linearly with the coarsening level.
 We  expect that according to Theorem \ref{Theo}ii, the error should be independent of the coarse graining parameter due to the microscopic nature of the algorithm though this is not evident in the tables since we are using a simplification of the reconstruction procedure for computational ease. 
We should also mention that a large number of samples ($10^5 $) were considered   ensuring   the statistical error is small enough.
 
\begin{table}[ht]
\caption{
N=  1028,
K= 1,
J= 5, 
 $\text{Error}_{cl}  = 0.003$,   Classical CPU = 94.5min
}
 \begin{tabular}{lcc} 
\hline\noalign{\smallskip}  
 CG level  & $\text{Error}_{c}$ &   Coupled CPU(min)     \\
 \noalign{\smallskip}\svhline\noalign{\smallskip}
q=4 &  0.01   &  23.1     \\
q=8 & 0.04    &    12.1   \\
\noalign{\smallskip}\hline\noalign{\smallskip}
 \end{tabular}
 \label{cost1}
 \centering
 \end{table}

%%%%%%%%%%%%%%%%%%%%%%%%%%%%%%%%%%%%%%%%%%%%%%%%%%%%%%%%%%%%%%%%%%%%%%%%%%%%%
%CONCLUSIONS
%%%%%%%%%%%%%%%%%%%%%%%%%%%%%%%%%%%%%%%%%%%%%%%%%%%%%%%%%%%%%%%%%%%%%%%%%%%%%
\section{Conclusions}

An advantage of the Coupled CGMC approach over the asymptotics methodology discussed in Section \ref{lattice-sl} is that the trial distribution  may even be order one away from the  target distribution, however, the method can  still perform well. On the other hand,      the methods can {\em complement} each other;  for example,  for  equilibrium sampling  considered in this work we  use   as a  trial reconstructed distribution,  the conditional measure $\nu(d\sigma|\eta)$ in  the multiscale decomposition in \cite{KPRT3}, see also    Section \ref{recon}.
 Such proposals based on careful statistical mechanics-based approximations    provide better trial choices for the MH methods and more efficient sampling, as is proved theroretically and numerically.
   The example illustrated makes clear that the coupled CGMC method implements a splitting of the short and long-range interaction terms, into the two Metropolis acceptance criteria involved. The long-range part which  is responsible for the expensive calculations at a fully  microscopic method, now enters only in the  coarse approximation measure where its  computational cost is much lower.
  
   Coupling  of a coarse and fine step is  also   effective in the study of dynamic processes of stochastic lattice systems  with kinetic Monte Carlo methods, a topic    studied in detail in \cite{KKP}.

%%%%%%%%%%%%%%%%%%%%%%%%%%%%%%%%%%%%%%%%%%%%%%%%%%%%%%%%%%%%%%%%%%%%%%%%%%%%%
%REFERENCES
%%%%%%%%%%%%%%%%%%%%%%%%%%%%%%%%%%%%%%%%%%%%%%%%%%%%%%%%%%%%%%%%%%%%%%%%%%%%%

\bibliographystyle{spmpsci}
\bibliography{reference}
 
\end{document}